\documentclass[12pt]{amsart}
\usepackage{amsmath}
\usepackage{amsfonts,amssymb,amsthm, txfonts, pxfonts,amscd}
\def\struckint{\mathop{%
\def\mathpalette##1##2{\mathchoice{##1\displaystyle##2}%
  {##1\textstyle##2}{##1\scriptstyle##2}{##1\scriptscriptstyle##2}}%
\mathpalette
{\vbox\bgroup\baselineskip0pt\lineskiplimit-1000pt\lineskip-1000pt
\halign\bgroup\hfill$}
{##$\hfill\cr{\intop}\cr\diagup\cr\egroup\egroup}%
}\limits}

\newcommand{\specrad}[1]{\mathcal{R}(#1)}

\newtheorem*{thma}{Theorem A}
\newtheorem*{thmb}{Theorem B}
\newtheorem{theorem}{Theorem}[section]
\newtheorem{lemma}[theorem]{Lemma}
\newtheorem{corollary}[theorem]{Corollary}

\newtheorem{definition}[theorem]{Definition}

\newtheorem{notation}[theorem]{Notation}
\newtheorem{theorem-definition}[theorem]{Theorem-Definition}

\theoremstyle{remark}

\newcommand{\cx}{\mathbb{C}}
\newcommand{\integers}{\mathbb{Z}}

\newcommand{\reals}{\mathbb{R}}

\DeclareMathOperator{\diag}{Diag}

\DeclareMathOperator{\tr}{tr}

\DeclareMathOperator{\Sp}{Sp}
\DeclareMathOperator{\SL}{SL}

\providecommand{\matnorm}[1]{\lvert\lVert#1\rVert\rvert}
\providecommand{\opnorm}[1]{\lvert\lVert#1\rVert\rvert_{\text{op}}}
\providecommand{\norm}[1]{\lVert#1\rVert}

\providecommand{\abs}[1]{\lvert#1\rvert}
\providecommand{\dotprod}[2]{\langle#1,#2\rangle}

\begin{document}
\title[Walks and bounds]{Walks on graphs and lattices -- effective bounds and applications}

\author{Igor Rivin}
\address{Department of Mathematics, Temple University, Philadelphia}
\email{rivin@math.temple.edu}
\thanks{The author would like to thank Peter Sarnak for encouragement
  and interesting conversations, and Hebrew University, University of
  Chicago, and Stanford  University for their hospitality during the
  preparation of this paper. Most of the results presented in this
  paper were circulated in a preprint in the Summer of 2006}

\curraddr{Mathematics Department, Stanford University, Stanford,
  California}
\email{rivin@math.stanford.edu}
\date{\today}

\begin{abstract}
We continue the investigations started in \cite{walks,rirred}. We
consider the following situation: $G$ is a finite directed graph,
where to each vertex of $G$ is assigned an element of a finite group
$\Gamma.$ We consider all walks of length $N$ on $G,$ starting from
$v_i$ and ending at $v_j.$ To each such
walk $w$ we assign the element of $\Gamma$ equal to the product of the
elements along the walk. The set of all walks of length $N$ from $v_i$
to $v_j$ thus induces a probability distribution $F_{N, i, j}$ on
$\Gamma.$ In \cite{walks} we give necessary and sufficient conditions
for the limit as $N$ goes to infinity of $F_{N, i, j}$ to exist and to
be the uniform density on $\Gamma$ (a detailed argument is presented
in \cite{rirred}). The convergence speed is then exponential in $N.$

In this paper we consider $(G, \Gamma),$ where $\Gamma$ is a group
possessing Kazhdan's property $T$ (or, less restrictively, property
$\tau$ with respect to representations with finite image), and a
family of homomorphisms $\psi_k: \Gamma \rightarrow \Gamma_k$ with
finite image. Each $F_{N, i, j}$ induces a distribution
 $F_{N, i, j}^k$ on $\Gamma_k$ (by push-forward under $\psi_k$). Our
 main result is that, under mild technical assumptions, the
 exponential rate of convergence of $F_{N, i, k}^k$ to the uniform
 distribution on $\Gamma_k$ does not depend on $k.$ 

As an application, we prove effective versions of the results of
\cite{rirred} on the probability that a random (in a suitable sence)
element of  $\SL(n, \integers)$ or $\Sp(n, \integers)$ has irreducible
characteristic polynomial, generic Galois group, etc.
\end{abstract}
\maketitle

\section*{Introduction}

The following set-up was first brought up in \cite{walks}, and then
fleshed out and applied in a somewhat unexpected direction in
\cite{rirred}:

Firstly, let $G$ be a finite ``ergodic'' undirected graph, which means
that the adjacency matrix of $G$ has a unique Perron-Frobenius
eigenvalue with a strictly positive eigenvector.

Secondly, let $\Gamma$ be a finite group, and assign to each vertex
$v$ of $G$ an element $\gamma(v) \in \Gamma.$

Finally, consider the set of walks $W_{N, i, j}$ on $G$ of length $N$
starting at $v_i$ and ending at $v_j.$ Each walk $w \in W_{N, i, j}$
defines an element $\gamma(w) \in \Gamma:$ the element $\gamma(w)$ is
simply the product (in order) of elements $\gamma(v)$ along $w.$ The
set $W_{N, i, j}$ thus induces a probability distribution 
$F_{N, i,  j}$ on $\Gamma,$ where the probability $p_{N, i, j}(\nu)$ assigned to $\nu \in
\Gamma$ is defined as:
\[
p_{N, i, j}(\nu) = \dfrac{|\{w\in W_{N, i, j} \left| \gamma(w) =
    \nu\right.\}|}{|W_{N, i, j}|}.
\]
\emph{A priori,} it is not clear that $F_{N, i, j}$ ever has full
support, but, rather surprisingly, the following holds:
\begin{thma}[\cite{walks,rirred}]
\label{mythm2}
If the set $\{\gamma(v)\left| v\in V(G)\right.\}$ generates $\Gamma$ and
  there is no one-dimensional complex representation $\rho$ of
  $\Gamma$ which maps all of $\gamma(v)$ to the same complex number,
  then the distributions $F_{N, i, j}$ converge to the \emph{uniform}
  distribution on $\Gamma.$ The speed of convergence is exponential in $N.$
\end{thma}
The proof of Theorem A is recalled below. The application to
irreducibility of random matrices in \cite{rirred} requires the use of
Theorem A for finite quotients of $\SL(n, \integers)$ and $\Sp(2n,
\integers).$ To get effective bounds, we need to have uniform bounds
on the exponential speed of convergence in Theorem A, and this is the
main subject of the current paper. The setup is as before, but $\Gamma$ is no longer
(necessarily) finite, but it \emph{is} assumed to have property $\tau$
for representations with finite image (see \cite{lubotztau} for
discussion of Property $\tau$). Any finite homomorphism $\psi$ of
$\Gamma$ with finite image $\psi{G}$ induces a family of distributions
$F_{N, i, j}^{\psi}$ on $\psi{G}.$ We then have the following:
\begin{thmb}
Let $G,\Gamma$ be as above. With the assumptions as in Theorem A, and
the additional assumption that the set $\{\gamma(v)^{-1} \gamma(w)
\left| v, w \in V(G)\right.\}$ generates $\Gamma$ the exponential
convergence rate of $F_{N, i, j}^\psi$ to the uniform distribution on
$\psi(\Gamma)$ can be bounded independently of $\psi.$
\end{thmb}

The plan of the rest of the paper is as follows:

The starting point for the proof of the theorems above is Fourier
Transform on finite groups, which is discussed in Section
\ref{fouriergroups}. In particular, we will be using Theorem
\ref{fourierest} and Corollary \ref{closetoconst} to reduce the
question of whether a probability distribution is close to uniform to
the proving that the Fourier Transform is small at every
\emph{non-trivial} representation. The reader might well wonder how
moving the problem to Fourier transform space helps us -- the answer
is that it turns out that we can reduce the estimation of the
``fourier coefficients'' to questions in linear algebra, through the
construction in Section \ref{fouriertomat}.

In Sections \ref{kazht}, \ref{linalgest}, \ref{shrinkageapp} we prove
the additional estimates we need to prove Theorem B. Finally, in
Section \ref{irredapp} we use Theorem B to show that the probability 
that a matrix in $\SL(n, \integers)$ or in $\Sp(2n, \integers)$ given
by a word of length $N$ in a symmetric generating set has reducible
characteristic polynomial decreases exponentially with $N.$

\section{Fourier Transform on finite groups}
\label{fouriergroups}
For a thorough introduction to the topic of this section the reader is
referred to \cite{serrereps,simonreps}.
Let $\Gamma$ be a finite group, and let $f:\Gamma\rightarrow \cx$ be a
function on $\Gamma.$ Furthermore, let $\widehat{\Gamma}$ be the
\emph{unitary dual} of $\Gamma:$ the set of all irreducible complex unitary
representations of $\Gamma.$ To $f$ we can associate its \emph{Fourier
  Transform} $\hat{f}.$ This is a function which associates to
each $d$-dimensional unitary representation $\rho$ a $d\times d$ matrix
$\hat{f}(\rho)$ as follows:
\[
\hat{f}(\rho) = \sum_{\gamma \in \Gamma} f(\gamma) \rho(\gamma).
\]

There is an inverse transformation, as well. Given a function $g$
on $\widehat{\Gamma}$ which associates to each $d$-dimensional
representation $\rho$ a $d\times d$ matrix $g(\rho),$ we can write:
\[
g^\sharp (\gamma) = \dfrac{1}{\abs{\Gamma}} \sum_{\rho\in
\widehat{\Gamma}} d_\rho \tr (g(\rho) \rho(\gamma^{-1}),
\]
where $d_\rho$ is the dimension of $\rho.$ We mean ``inverse'' in the
most direct way possible:
\[
\hat{f}^\sharp = f.
\]

The following result is classical (see, eg, \cite{simonreps}):
\begin{theorem}
\[
\sum_{\rho in \widehat{\Gamma}} d_\rho^2 = \abs{\Gamma},
\]
\end{theorem}
and, together with the Fourier inversion formula, implies
\begin{theorem}
\label{fourierest}
Let $g$ be a function on $\widehat{\Gamma},$ such that for every
\emph{nontrivial} $\rho \in \widehat{\Gamma},$ 
\[
\opnorm{g(\rho)} < \epsilon,
\]
where $\opnorm{\bullet}$ denotes the operator norm (see Section
\ref{matnorm}). 
Then, for any $\gamma_1, \gamma_2 \in \Gamma,$
\[
\abs{g^\sharp(\gamma_1) - g^\sharp(\gamma_2)}  < 2\epsilon.
\]
\end{theorem}
\begin{proof}
First, note that for the trivial representation $\rho_0,$ the quantity
\[
d_{\rho_0}g(\rho_0)\rho_0(\gamma) = g(\rho_0),
\]
so does not depend on $\gamma.$
By the Fourier inversion formula, then,
\[
\begin{split}
\abs{g^\sharp(\gamma_1) - g^\sharp(\gamma_2)} & = \\
\left\lvert \dfrac{1}{\abs{\Gamma}}\sum_{\substack{\rho \in \widehat{\Gamma}\\
       \rho \neq \rho_0}} d_\rho \tr(g(\rho) (\rho(\gamma_1) -
   \rho(\gamma_2)))\right\rvert & \leq \\
\sum_{i=1}^2 \left\lvert \dfrac{1}{\abs{\Gamma}}\sum_{\substack{\rho \in \widehat{\Gamma}\\
       \rho \neq \rho_0}} d_\rho \tr\left(g(\rho) \rho(\gamma_i) \right)
   \right\rvert &
\underset{\text{by
     Eq. \eqref{traceineq}}}{\leq} \\
\dfrac{2}{\abs{\Gamma}} \sum_{\rho \in \widehat{\Gamma}}  d_\rho^2
\opnorm{g(\rho)} & < 2\epsilon.
\end{split}
\]
\end{proof}
\begin{corollary}
\label{closetoconst}
Under the assumption of Theorem \ref{fourierest}, and assuming in
addition that $g$ is real valued, if 
\[\sum_{\gamma \in \Gamma} g(\gamma) = 1,
\]
then \[g(\gamma) - 1/|\Gamma| < 2 \epsilon\quad\forall \gamma \in
\Gamma.\]
Furthermore, if $\Omega \in \Gamma,$ 
\begin{equation}
\label{closetoconstset}
\left| \sum_{\gamma \in \Omega} g(\gamma) -
  \dfrac{\Omega}{|\Gamma|}\right| < 2 \epsilon |\Omega.|
\end{equation}
\end{corollary}
\begin{proof}
Without loss of generality, suppose that $g(\gamma) > 1/|\Gamma|.$
Then there is a $\gamma_2,$ such that $g(\gamma_2) < 1/|\Gamma|.$
Thus,
\[
g(\gamma) - 1/|\Gamma| < g(\gamma)-g(\gamma_2) < 2\epsilon.
\]
The estimate \eqref{closetoconstset} follows immediately by summing
over $\Omega.$
\end{proof}

\section{Fourier estimates via linear algebra}
\label{fouriertomat}
In order to prove Theorem A, we would like to use Theorem
\ref{fourierest}, and to show the equidistribution result, we would
need to show that for every \emph{nontrivial} irreducible
representation $\rho,$ 
\begin{equation}
\label{coeffdecay}
\lim_{N\rightarrow \infty}\dfrac{1}{|W_{N, i, j}|}\tr{\sum_{w\in W_{N, i,
    j}}  \rho(\gamma_w)} = 0.
\end{equation}

To demonstrate Eq. \eqref{coeffdecay}, suppose that $\rho$ is
$k$-dimensional, so acts on a $k$-dimensional Hilbert space $H_\rho=H.$
Let $Z = L^2(G)$ -- the space of complex-valued functions from $V(G)$
to $\cx,$ let $e_1, \dotsc, e_n$ be the standard basis of $Z,$ and let
$P_i$ be the orthogonal projection on the $i$-th coordinate space. We
introduce the matrix 
\[
U_\rho = \sum_{i=1}^n P_i \otimes \rho(t_i) = 
\begin{pmatrix}
\rho(t_1) & 0 & \dots & 0\\
0 & \rho(t_2) & \dots & 0\\
\hdotsfor[2]{4}\\
0 & 0 & \dots & \rho(t_n)
\end{pmatrix},
\]
and also the matrix $A_\rho = A(G) \otimes I_H,$ where $I_H$ is the
identity operator on $H.$ Both $U_\rho$ and $A_\rho$ act on $Z \otimes H.$
The following is immediate:
\begin{lemma}
\label{matprodlem}
Consider the matrix $(U_\rho A_\rho)^l,$ and think of it as an
$n\times n$ matrix of $k\times k$ blocks. Then the $ij$-th block
equals the sum over all paths $w$ of length $l$ beginning at $v_i$ and
ending of $v_j$ of $\rho(\gamma_w).$
\end{lemma}

Now, let $T_{ji}$ be the operator on $Z$ which maps $e_k$ to
$\delta_{kj} e_i.$

\begin{lemma}
\label{projlemma}
\[
\tr{\left[\left((T_{ji}^t P_j)\otimes I_H\right) (U_\rho \otimes A_\rho)^N (P_i\otimes I_H)\right]} = \tr{\sum_{w\in  W_{N, i, j}} \rho(\gamma_w)}
\]
\end{lemma}

\begin{proof}
The argument of trace on the left hand side simply extracts the $ij$-th $k\times k$ block from $(U_{\rho} \otimes A_\rho)^N.$
\end{proof}

By submulticativity of operator norm, we see that 
\[
\opnorm{(T_{ji}^t P_j)\otimes I_H (U_\rho \otimes A_\rho)^N P_i\otimes
  I_H} \leq \opnorm{ (U_\rho \otimes A_\rho)^N}, 
\]
and so proving Theorem A reduces (thanks to Theorem
\ref{fourierest}) to showing
\begin{theorem}
\label{fundcollapse}
\[
\lim_{N\rightarrow \infty} \dfrac{\opnorm{(U_\rho \otimes
    A_\rho)^N}}{|W_{N, i, j}|
%\opnorm{A_\rho^N}
} = 0,
\]
for any non-trivial $\rho.$
\end{theorem}
\begin{notation}
We will denote the spectral radius of an operator $A$ by $\specrad{A}.$
\end{notation}
Since $|W_{N, I, j}| \asymp \mathcal{R}^N(A(G)),$ and by Gelfand's Theorem
(Theorem \ref{gelfand}), 
\[
\lim_{N\rightarrow \infty} \|B^N\|^{1/N} = \specrad{B},
\]
for any matrix $B$ and any matrix norm $\|\bullet\|,$ Theorem
\ref{fundcollapse} is equivalent to the statement that the spectral
radius of $U_\rho \otimes A_\rho$ is smaller than that of $A(G).$

Theorem \ref{fundcollapse} is proved in Section \ref{myproof}. 

\subsection{Proof of  Theorem  \ref{fundcollapse}}
\label{myproof}
\begin{lemma}
\label{twistlemma}
Let $A$ be a bounded hermitian operator $A:H\rightarrow H,$ and $U:
H\rightarrow H$ a unitary operator on the same Hilbert space $H.$
Then the spectral radius of $U A$ is smaller than the spectral radius
of $A,$ and the inequality is strict unless an eigenvector of $A$ with
maximal eigenvalue is also an eigenvector of $U.$
\end{lemma}
\begin{proof}
The spectral radius of $UA$ does not exceed the operator norm of $UA,$
which is equal to the spectral radius of $A.$ Suppose that the two are
equal, so that there is a $v,$ such that $\norm{UA v} = \specrad{A}{v},$
and $v$ is an eigenvector of $UA.$
Since $U$ is unitary, $v$ must be an eigenvector of $A,$ and since it
is also an eigenvector of $UA,$ it must also be an eigenvector of $U.$
\end{proof}

In the case of interest to us, $\rho$ is a $k$-dimensional irreducible
representation of $\Gamma,$ $U = \diag(\rho(t_1), \dots,
\rho(t_n),$ while $A = A(G) \otimes I_k.$ We assume that $A(G)$ is an
irreducible matrix, so that there is a unique eigenvalue of modulus
$\specrad{A(G)},$ that eigenvalue $\lambda_{\max}$ 
(the \emph{Perron-Frobenius eigenvalue})
is positive, and it has a strictly positive eigenvector $v_{\max}.$ We
know that the spectral radius of $A$ equals the spectral radius of
$A(G),$ and the eigenspace of $\lambda_{\max}$ is the set of vectors of
the form $v_{\max} \otimes w,$ where $w$ is an arbitrary vector in
$\mathbb{C}^k.$ If $v_{\max} = (x_1, \dotsc, x_n),$ we can write
$v_{\max} \otimes w = (x_1 w, \dotsc, x_n w),$ and so $U(v_{\max} \otimes
w) = (x_1 \rho(t_1) w, \dotsc, x_n \rho(t_n) w).$ Since all of the
$x_i$ are nonzero, in order for the inequality in Lemma
\ref{twistlemma} to be nonstrict, we must have some $w$ for which
$\rho(t_i) w = c w$ (where the constant $c$ does \emph{not} depend on
$i.$) 
Since the elements $t_i$ generate $\Gamma,$ the existence of such a
$w$ contradicts the irreducibility of $\rho,$ \emph{unless} $\rho$ is
one dimensional. This proves Theorem \ref{fundcollapse}

\section{Some remarks on matrix norms}
\label{matnorm}
In this note we use a number of matrix norms, and it is useful to
summarize what they are, and some basic relationships and inequalities
satisfied by them. For an extensive discussion the reader is referred
to the classic \cite{hornjohnson}. All matrices are assumed square,
and $n\times n.$

A basic tool in the inequalities below is the \emph{singular value
decomposition} of a matrix $A.$ 
\begin{definition}
The singular values of $A$ are the non-negative
square roots of the eigenvalues of
$A A^*,$ where $A^*$ is the conjugate transpose of $A.$ 
\end{definition}
Since $AA^*$
is a positive semi-definite Hermitian matrix for any $A,$ the singular
values $\sigma_1 \overset{\text{def}}{=} \sigma_{\max} \geq \sigma_2
\geq \dots$ are non-negative real numbers. For a Hermitian $A,$ the
singular values are simply the absolute values of the eigenvalues of $A.$

The first matrix norm is the \emph{Frobenius norm}, denoted by
$\norm{\bullet}.$
This is defined as 
\[
\norm{A} = \sqrt{\tr{A A^*}} = \sqrt{\sum_i \sigma_i^2}.
\]
This is also the sum of the square moduli of the elements of $A.$

The next matrix norm is the \emph{operator norm}, $\opnorm{\bullet},$
defined as 
\[
\opnorm{A} = \max_{\norm{v} = 1} \norm{Av} = \sigma_{\max}
\]

Both the norms $\norm{\bullet}$ and $\opnorm{\bullet}$ are
\emph{submultiplicative} (submultiplicativity is part of the
definition of matrix norm: saying that the norm $\matnorm{\bullet}$ is
submultiplicative means that $\matnorm{AB} \leq \matnorm{A}
\matnorm{B}$.)

From the singular value interpretation\footnote{A celebrated result of
  John von Neumann states that \emph{any} unitarily invariant matrix
  norm is a symmetric guage on the space of singular values -
  \cite{vnguage}.} of the two matrix norms and the Cauchy-Schwartz
inequality we see immediately that 
\begin{equation}
\label{twonorms}
\norm{A}/\sqrt{n}\leq \opnorm{A}\leq \norm{A}
\end{equation}

We will also need the following simple inequalities:
\begin{lemma}
\label{ulemma}
Let $U$ be a unitary matrix:
\begin{equation}
\label{traceineq}
\abs{\tr AU} \leq \norm{A}\sqrt{n} \leq n \opnorm{A}.
\end{equation}
\end{lemma}
\begin{proof}
Since $U$ is unitary, $\norm{U}= \norm{U^t} = \sqrt{n}.$
So, by the Cauchy-Schwartz inequality, 
$tr A U \leq \norm{A} \norm{U} = \sqrt{n} \norm{U}.$
The second inequality follows from the inequality \eqref{twonorms}.
\end{proof}

The final (and deepest result) we will have the opportunity to use is:
\begin{theorem}[Gelfand]
\label{gelfand}
For any 
operator $M,$ the spectral radius $\specrad{M}$ and any matrix norm
$\matnorm{\bullet},$
\[
\specrad{M} = \lim_{k\rightarrow \infty} \matnorm{M^k}^{1/k},
\]
\end{theorem}

\section{Some remarks on Kazhdan's property T}
\label{kazht}
A group $G$ is said to have \emph{Kazhdan's Property $T$} if there
exists an $\epsilon>0$ and a compact subset $K \subseteq G$ such that
for every nontrivial irreducible representation $(H, \rho)$ of $G$ and
every vector $v\in H$ of norm one, $\|\rho(k) v - v\| > \epsilon$ for
some $k \in K.$ This definition is the one given in A. Lubotzky's book
\cite{lubotzbook}. For finitely generated discrete groups $K$ can be
taken to be any set of generators (though the
$\epsilon$\footnote{known as Kazhdan's constant} will depend
on the generating set, it is obvious that knowing Kazhdan's constant
for some generating set will give bounds for any other generating set.
It is known that lattices in semi-simple Lie groups have property $T$
and Kazhdan's constants have been explicitely computed by Y. Shalom
(see \cite{shalomkazhdan}). Related results have also been obtained by
A. Zuk \cite{zukkazhdan}.

We will need the following
\begin{lemma}
\label{tprime}
Let $G$ have Kazhdan's property $T$ and let $t_1, \dots, t_n$ be a
generating set of $G,$ such that the set of all products $t_j^{-1}t_i$
is also a generating set. Then, there exists an $\epsilon > 0$ such
that for any irreducible representation $(H, \rho)$ and any pair $v, w
\in H$ there exists $i\leq n$ such that $\|\rho(t_i) v - 
w\| > \epsilon.$
\end{lemma}
\begin{proof}
Suppose not.
By the triangle inequality, $\|\rho(t_i)v - \rho(t_j)v \| <
2\epsilon,$ for all pairs $i, j.$ Since $\rho$ is unitary, we see that 
$\|\rho(t_j^{-1}t_i) v - v\| < 2\epsilon.$ It follows that the we can
choose the $\epsilon$ whose existence is postulated in the Lemma to be
half the Kazhdan constant of $G$ with respect to the generating set
consisting of all products $t_j^{-1} t_i.$
\end{proof}

To show that the condition in the statement of Lemma \ref{tprime} is
often met, first note:
\begin{lemma}
\label{indextwo}
Let $S = \{t_1, \dotsc, t_n\}$ be a \emph{symmetric} generating set for $G.$
Then, the subgroup $H$ generated by all products of the form $t_j^{-1}
t_i$ has index at most two in $G$ (hence is always normal).
\end{lemma}
\begin{proof}
Since $S$ is symmetric, $H$ has every element which can be written as
a word of even length in the elements of $S.$ If $H \neq G,$ then the
index of $H$ clearly equal to two (the other coset being the set of
``odd'' elements of $G.$
\end{proof}
\begin{corollary}
If $G$ is one of \[
\SL(n, \mathbb{Z}), \SL(n, \mathbb{Z}/p\mathbb{Z}),
\Sp(n, \mathbb{Z}), \Sp(n, \mathbb{Z}/p \mathbb{Z})\] for $n \geq 2,$ 
and $S$ is a symmetric generating set, then $S^{-1}S$ generates $G.$
\end{corollary}

\section{Linear algebra estimates}
\label{linalgest}
\begin{lemma}
\label{twistang}
Let $U, A$ be as in Lemma \ref{twistlemma}. Assume that the spectral
radius of $A$ equals $1$ (for simplicity of notation), that the second
biggest (in absolute value) eigenvalue of $A$ has absolute value
$\lambda < 1.$ Let $A_{\max}$ be the eigenspace of $A$ corresponding
to the eigenvalue $1,$ and let $P_{\max}$ be the orthogonal projection
on $A_{\max}.$ Assume now that for any  $v\in A_{\max},$ 
\begin{equation}
\label{distineq}
\norm{P_{\max} U v} \leq d \norm{v},
\end{equation}
for some $0\leq d < 1.$
Then, there is a function $f(\lambda, d)<1,$ such that the spectral
radius of $UA$ is smaller than $f(\lambda, d).$
\end{lemma}

\begin{proof}
We will use Gelfand's Theorem \ref{gelfand}
For our result, we will use the operator norm, and Lemma
\ref{twistang} will follow immediately from Theorem \ref{twistshrink},
with $f(\lambda, d) = \sqrt{g(\lambda, d)},$ where $g$ is the function
in the statement of Theorem \ref{twistshrink}.
\end{proof}

\begin{theorem}
\label{twistshrink}
For $U, A$ as in the statement of Lemma \ref{twistang}, and $v$ an
arbitrary vector. Then \[\norm{(UA)^2v} \leq g(\lambda, d)\norm{v},\] for some
function $g(\lambda, d) < 1,$ and so \[\opnorm{(UA)^k} \leq 
g^{\lfloor  k/2 \rfloor}(\lambda,d),\] where $\opnorm{M}$ denotes the operator
norm of $M.$
\end{theorem}

\begin{proof}
Since $U$ is unitary, $\norm{(UA)^2v} = \norm{AUAv},$ for any $v.$
Now write $v=x \oplus y,$ with $x \in A_{\max},$ and $y \in
A_{\max}^\perp.$

Our first observation is that 
\begin{equation}
\label{lenchange}
\norm{Av}^2 \leq \norm{x}^2 + \lambda^2 \norm{y}^2 = \lambda^2 \norm{v}^2 +
(1-\lambda^2) \norm{x}^2.
\end{equation}
It follows that 
\begin{equation}
\label{bigineq1}
\norm{A U A} \leq \norm{A}.
\end{equation}

Our second observation is that 
\begin{equation}
\label{unichange}
\norm{P_{\max} U A v} \leq d \norm{x} + \lambda \norm{y},
\end{equation}
and so by \eqref{lenchange},
\begin{equation}
\label{allchange}
\begin{split}
\norm{A U A v}^2   & \leq \\
\lambda^2 \norm{A v}^2 + (1-\lambda^2) (d \norm{x} +
\lambda \norm{y})^2 & \leq \\
\lambda^2 (\norm{x}^2 + \lambda^2\norm{y}^2) +
(1-\lambda^2) (d \norm{x} + \lambda \norm{y})^2  & = \\
(1-(1-d^2)(1-\lambda^2) ) \norm{x}^2 + \lambda^2 \norm{y}^2 +
2(1-\lambda^2) d \lambda y x.
\end{split}
\end{equation}

Let us now write $\norm{y} = \alpha \norm{x}.$ 

\medskip\noindent\textbf{$\lambda > 0.$}

Eq. \eqref{allchange} gives
us 
\begin{equation}
\begin{split}
\dfrac{\norm{A U A v}^2}{\norm{v}^2} & = \\
\dfrac{1-(1-\lambda^2)(1-d^2) + \lambda^2 \alpha^2 +
  2(1-\lambda^2) d \lambda \alpha}{1+\alpha^2} & \leq \\
1- (1-\lambda^2) (1-d^2) + \lambda^2 \alpha^2 +
  2(1-\lambda^2) d \lambda \alpha = h(\lambda, d, \alpha).
\end{split}
\end{equation}
Note that $h(\lambda, d, 0) = 1-(1-\lambda^2)(1-d^2) < 1,$
and $h(\lambda, d, \alpha)$ is a monotonically increasing function of
$\alpha$ when $\alpha \geq 0,$ and $0\leq \lambda, d< 1.$
This means that we can find $0<\alpha_0$ such that 
$h(\lambda, d, \alpha_0) = 1-(1-\lambda^2)(1-d^2)/2,$
namely 
\begin{equation}
\label{alphadef}
\alpha_0 = \dfrac{1-\lambda^2}{\lambda} \left(\sqrt{d^2 +
    \dfrac{1-d^2}{2(1-\lambda^2)}}\right),
\end{equation}
Putting together all the inequalities, we see that if 
$\norm{y}/\norm{x} \leq \alpha_0,$ then 
\[
\norm{UAUAv} \leq \sqrt{1-(1-\lambda^2)(1-d^2)/2} \norm{v},
\]
while if $\norm{y}/\norm{x} > \alpha_0,$ then
\[
\norm{UAUAv} \leq \sqrt{\dfrac{1+\alpha_0 \lambda}{1+\alpha_0}} \norm{v},
\]
so setting
\[
g(\lambda, d) = \min\left(\sqrt{\dfrac{1+\alpha_0
      \lambda}{1+\alpha_0}}, 
  \sqrt{\dfrac{1-(1-\lambda^2)(1-d^2)/2}{1+\alpha_0}}\right),
\]
 the Lemma is proved.

\medskip\noindent\textbf{$\lambda=0.$}
In this case, the computation is much simpler:
\begin{equation}
\dfrac{\norm{A U A v}^2}{\norm{v}^2} = \dfrac{d^2}{1+\alpha^2}\leq d^2,
\end{equation}
and so the Lemma is proved here too.
\end{proof}

\section{Applications of Theorem \ref{twistshrink} to speed of
  convergence in Theorem A}
\label{shrinkageapp}
Let us apply Theorem \ref{twistshrink} to the setting of Theorems
A and B. We will be using the argument and the notation of 
Sections \ref{myproof} and \ref{linalgest}.
Let $S=\{t_1, \dotsc, t_n\},$ let $\Gamma$ be the group generated by
$S,$ and let $\Gamma_1$ be the group generated by $S^{-1} S.$

If $\lambda_1$ is the Perron-Frobenius eigenvalue of $G,$ and
$\lambda_2$ is the second largest (in absolute value) eigenvalue, we
set $\lambda=|\lambda_2|/|\lambda_1|.$
Let $X=(x_1, \dotsc, x_n)$ be the (unit) Perron-Frobenius eigvenctor of $A(G).$
We know that $A_1$ is the space of all vectors of the form
$Y=X\otimes v = (x_1 v, \dotsc, x_n v),$ where $v\in \reals^k.$
Such a vector is a unit vector precisely if $\norm{v} = 1.$
Recall that $U Y = (x_1 \rho(t_1) v, \dots, x_n \rho(t_n) v).$ Let 
$W= X \otimes w \in A_1,$ then 
\begin{equation}
\label{dotprodeq}
\dotprod{UY}{W} = \sum_{i=1}^n x_i^2 \dotprod{\rho(t_i) v}{w}.
\end{equation}

Assume now that the group $\Gamma_1$ has the analogue of 
Kazhdan's property $T,$ but with
respect to the set of restrictions of irreducible representations of
$\Gamma$ -- these are not necessarily irreducible when restricted to
$\Gamma_1$ --  with
the constant $\epsilon_1$ corresponding to the generating set $S^{-1}S.$ 
We know (by Lemma \ref{tprime}) that there is an $i\leq n,$ such that
$\|\rho(t_i) v - w\| \geq \epsilon_1/2,$ and so, by the Law of
Cosines, 
\[\dotprod{\rho(t_i) v}{w} \leq 1 - \epsilon_1^2/8, \]
and so, by Eq. \eqref{dotprodeq}, 
\[
\dotprod{UV}{W} \leq 1 - x_i^2\epsilon_1^2/8
\]
Lemma \ref{twistang} now gives us:
\begin{lemma}
\label{myeff}
The operator norm  of $(UA)^k$ is at most $g^{\lfloor k/2\rfloor}(\lambda, 1- x_i^2
\epsilon_1^2/8),$ where $g$ is the function computed in Theorem \ref{twistshrink}.
\end{lemma}

This completes the proof of Theorem B.

\section{Applications to irreducibility}
\label{irredapp}
In this section, Theorem B is used to
show that the probability that a random walk of length $N$ on a graph
$G$ decorated with elements of $\SL(n, \integers)$ or $\Sp(2n,
\integers)$ represents a matrix with \emph{reducible} characteristic
polynomial goes to $0$ exponentially fast with the length $N$ of the
walks considered.

The results above show that for a fixed graph $G$ and the series of groups
$\Gamma_p,$ where $\Gamma_l = \SL(n, ;)$ or $\Gamma_l = \Sp(2n, l)$
there exist a constant $c>1,$ such that the probability $p_\gamma$ that one of the
random walks of length $N$ over $G$ (decorated with elements of
$\Gamma_p$) hits a subset $\Omega \subseteq \Gamma_l$  satisfies 
\begin{equation}
\label{omegaest}
|p_\Omega - |\Omega|/|\Gamma_l| \leq 2 c^{-N} |\Omega|,
\end{equation}
 where $c>1$ does \emph{not} depend on $l.$ 

\subsection{$\SL(n).$}
\label{sleff}
We know (see \cite{rirred}) that the set $\mathcal{R}_l \in \SL(n, l)$ has cardinality
bounded by 
\begin{equation}
\label{slest}
|\mathcal{R}_p| \leq c_2 |\SL(n, p)|,
\end{equation}
for $p$ prime. Now, for given $N\gg 1,$ there is a prime $p_N$ satisfying 
\[
(1-\epsilon)c^{N/(n^2-1)} \leq p_N \leq (1+\epsilon)c^{N/(n^2-1)}.
\]
By estimates \eqref{omegaest} and \eqref{slest}, it follows that a
random walk on $G$ of length $N$ represents a reducible element in
$\SL(n, p_N)$ with probability $P_N$ bounded above by:
\begin{equation}
\label{probslest}
P_N \leq \dfrac{c_2}{p_N}(1+(1+\epsilon)c_2) = O(c^{N/(n^2-1)}).
\end{equation}
Since an element in $\SL(n, \integers)$ is reducible only if it is
reducible in $\SL(n, l)$ (for every $l$), \eqref{probslest} gives an
upper bound on the probablity that an element represented by a random
walk of length $N$ is reducible over the integers.

\subsection{$\Sp(2n)$}
\label{speff}
Here, the method in the last section does not work (since we only have
$O(1)$ bounds for individual primes).

Therefore, define
\[
q_k = \prod_{i=1}^k p_k
\] (so $q_k$ is the product of the first $k$ primes).
The prime number theorem tells us that $q_k \sim k^k.$\footnote{If we
  wished to keep this discussion completely elementary, Chebyshev's
  elementary bound tells us that $q_k = O(k^{ak})$ for some $a > 1,$
  which is sufficient for what we are about to do.}

By Borel's estimate and the strong approximation property for
$\Sp(2n)$ (see \cite{rirred}) we
know that probability that an element of $\Sp(2n, q_k)$ is reducible
is bounded above by $c_3^k,$ for some $c_3 < ,$ and so by
\eqref{omegaest} we know that the probability $P_N$ that a walk on $G$
of length $N$ gives us a reducible element modulo $q_k$ is bounded above by
\[
P_N \leq c_3^{-k}(1+ 2 c^{-N} k^{k (2n^2+n)}).
\]
If we pick 
\[k \approx \dfrac{\dfrac{N}{2n^2+n}\log c }{\log\dfrac{N}{2n^2+n}}\]
(so that the second term in parenthesese is $O(1)$), we see that 
\[
P_N = O(\exp(\log c_3 \log c (N/(2n^2+n) - \epsilon))),
\]
for any $\epsilon > 0,$ and as before, the same bound obtains for the
probability that a random walk of length $N$ on $G$ gives a reducible
element in $\Sp(2n, \integers).$

\subsection{Remarks}

The first observation is that the argument in Section \ref{speff}
applies, \emph{mutatis mutandis} to the problem of counting elements
in $\Sp(n, \integers)$ whose Galois group is not the full symmetric
group.

Secondly, presumably sharper bounds can be given using more
sophisticated sieve machinery (see, eg, \cite{MurtySieve}). As
evidence for this, if the argument above is used to estimate the
probability that a 
\emph{polynomial} of degree $d$ with coefficient height bound $H,$
reducible, our argument gives $O(H^{\log(d-1) - \log d}),$ 
Gallagher's large sieve argument \cite{galgal} gives $O(H^{-1/2}),$
while the truth is $O(1/H).$ Since the arguments above are completely
elementary (even the use of the Prime Number Theorem can be avoided),
and we get the result we want (that the probability decays
exponentially) it seems wise to leave sieve methods to the experts.
In fact, related results have been obtained by Emmanuel Kowalski, using his
deep generalization of the large sieve \cite{kowalskisieve} (also
monograph, in preparation).
\bibliographystyle{plain}
\bibliography{rivin}
\end{document}